\documentclass[11pt]{amsart}

\usepackage{amscd,latexsym,amsthm,amsfonts,amssymb,amsmath,amsxtra}
\usepackage[mathscr]{eucal}
\usepackage[pagebackref,colorlinks=true,linkcolor=blue,citecolor=blue]{hyperref}
\usepackage{xcolor}

\renewcommand\theequation{\thesection.\arabic{equation}}

\newcommand{\BC}{{\mathbb {C}}}

\newcommand{\CH}{{\mathcal {H}}}

\newcommand{\CN}{{\mathcal {N}}}
\newcommand{\CO}{{\mathcal {O}}}

\newcommand{\CS}{{\mathcal {S}}}

\newcommand{\RG}{{\mathrm {G}}}

\newcommand{\GL}{{\mathrm{GL}}}

\newcommand{\G}{{\mathrm{G}}}

\newcommand{\SL}{{\mathrm{SL}}}

\newcommand{\SO}{{\mathrm{SO}}}

\newcommand{\Sp}{{\mathrm{Sp}}}

\newcommand{\tr}{{\mathrm{tr}}}

\newcommand{\wt}{\widetilde}

\newcommand{\ul}{\underline}

\newtheorem{thm}{Theorem}[section]

\newtheorem{lem}[thm]{Lemma}
\newtheorem{prop}[thm]{Proposition}
\newtheorem {conj}[thm]{Conjecture}

\newtheorem {ques/conj}[thm]{Question/Conjecture}

\newtheorem{rmk}[thm]{Remark}

\makeatletter

\newcommand{\Rmnum}[1]{\expandafter\@slowromancap\romannumeral #1@}
\makeatother

\begin{document}
\renewcommand{\theequation}{\arabic{equation}}
\numberwithin{equation}{section}

\title[Jiang conjecture]{The Jiang conjecture on the wavefront sets of local Arthur packets}

\author{Baiying Liu}
\address{Department of Mathematics\\
Purdue University\\
West Lafayette, IN, 47907, USA}
\email{liu2053@purdue.edu}

\author{Freydoon Shahidi}
\address{Department of Mathematics\\
Purdue University\\
West Lafayette, IN, 47907, USA}
\email{freydoon.shahidi.1@purdue.edu}


\date{\today}

\dedicatory{To Steve Kudla on the occasion of his 70th Birthday}


\thanks{The research of the first named author is partially supported by the NSF Grants DMS-1702218, DMS-1848058. The research of the second named author is partially supported by the NSF Grants DMS-1801273, DMS-2135021}

\begin{abstract}
This is a report on the progress made on a conjecture of Jiang on the upper bound nilpotent orbits in the wave front sets of representations in local Arthur packets of classical groups, which is a natural generalization of the Shahidi conjecture. We partially prove this conjecture, confirming the relation between the structure of wave front sets and the local Arthur parameters. Under certain assumptions, we also prove the enhanced Shahidi conjecture, which states that local Arthur packets are tempered if and only if they have generic members.  
\end{abstract}

\maketitle


\section{Introduction}

Let $F$ be a non-Archimedean local field. Let $\mathrm{G}_n=\Sp_{2n}, \SO_{2n+1}, \SO_{2n}^{\alpha}$ be  quasi-split classical groups, where $\alpha$ is a square class in $F$, and let $G_n=\mathrm{G}_n(F)$.  
Here, we identify a square class with the corresponding quadratic character of the Weil group $W_F$ via the local class field theory. 
Then the Langlands dual groups are 
$$\widehat{\mathrm{G}}_n(\BC) = \SO_{2n+1}(\BC), \Sp_{2n}(\BC), \SO_{2n}(\BC).$$
Let ${}^L\mathrm{G}_n$ be the $L$-group of $\mathrm{G}$,
$${}^L\mathrm{G}_n= 
\begin{cases}
\widehat{\mathrm{G}}_n(\BC) & \text{ when } \mathrm{G}_n=\Sp_{2n}, \SO_{2n+1},\\
\SO_{2n}(\BC) \rtimes W_F & \text{ when } \mathrm{G}_n=\SO_{2n}^{\alpha}.
\end{cases}
$$ 
In his fundamental work \cite{Art13}, Arthur introduces the local Arthur packets which are finite sets of representations of $G_n$, parameterized by local Arthur parameters. Local Arthur parameters are defined as
a direct sum of irreducible representations
$$\psi: W_F \times \SL_2(\mathbb{C}) \times \SL_2(\mathbb{C}) \rightarrow {}^L\mathrm{G}_n$$
\begin{equation}\label{lap}
  \psi = \bigoplus_{i=1}^r \phi_i \otimes S_{m_i} \otimes S_{n_i},  
\end{equation}
satisfying the following conditions:

(1) $\phi_i(W_F)$ is bounded and consists of semi-simple elements, and $\dim(\phi_i)=k_i$;

(2) the restrictions of $\psi$ to the two copies of $\SL_2(\mathbb{C})$ are analytic, $S_k$ is the $k$-dimensional irreducible representation of $\SL_2(\mathbb{C})$, and 
$$\sum_{i=1}^r k_im_in_i = N=N_n:= 
\begin{cases}
2n+1 & \text{ when } \mathrm{G}_n=\Sp_{2n},\\
2n & \text{ when } \mathrm{G}_n=\SO_{2n+1}, \SO_{2n}^{\alpha}.
\end{cases}
$$ 
Assuming the Ramanujan conjecture, Arthur (\cite{Art13}) showed that these local Arthur packets characterize the local components of square-integrable automorphic representations.
For $1 \leq i \leq r$, let $a_i=k_im_i$, $b_i=n_i$. 
Let
$$\ul{p}(\psi) = [b_1^{a_1} b_2^{a_2} \cdots b_r^{a_r}]$$
be a partition of $N$, where without loss of generality, we assume that $b_1 \geq b_2 \geq \cdots \geq b_r$. 

For each local Arthur parameter $\psi$, Arthur associated a local $L$-parameter $\phi_{\psi}$ as follows
\begin{equation}\label{apequ1}
\phi_{\psi}(w, x) = \psi\left(w, x, \begin{pmatrix}
        |w|^{\frac{1}{2}} & 0 \\
        0 & |w|^{-\frac{1}{2}}\\
\end{pmatrix}\right).
\end{equation}
Note that for any local Arthur parameter $\phi \otimes S_m \otimes S_n$,
$$\phi(w) \otimes S_m (x) \otimes S_n \left(\begin{pmatrix}
        |w|^{\frac{1}{2}} & 0 \\
        0 & |w|^{-\frac{1}{2}}\\
\end{pmatrix}
\right) = \bigoplus_{j=-\frac{n-1}{2}}^{\frac{n-1}{2}} |w|^j\phi(w) \otimes S_m (x)
.$$
Arthur also showed that $\psi \mapsto \phi_{\psi}$ is injective.
Let $\pi_{\psi}$ be the representation of $\GL_{N}(F)$ corresponding to $\phi_{\psi}$ via local Langlands correspondence, which is unitary and self-dual, and let $\wt{\pi}_{\psi}$ be its canonical extension to the bitorsor $\wt{\GL}_{N}(F) = {\GL}_{N}(F) \rtimes \wt{\theta}(N)$ 
where 
$\tilde{\theta}(N)=Int(\tilde{J})\circ \theta: g \mapsto \tilde{J} \theta(g) \tilde{J}^{-1}$, 
$\theta(g)={}^t g^{-1}$, 
$$
\tilde{J}=\tilde{J}(N)=
\begin{pmatrix}
0&&&1\\
&&-1&\\
&\dots&&\\
(-1)^{N+1}&&&0
\end{pmatrix},
$$
(see \cite[Section 2.2]{Art13} for more details).

Given an irreducible representation $\pi$ of $G_n$, one important invariant is a set $\frak{n}(\pi)$ which is defined to be all the $F$-rational nilpotent orbits $\CO$ in the Lie algebra $\frak{g}_n$ of $G_n$ such that the coefficient $c_{\CO}(\pi)$ in the Harish-Chandra-Howe local expansion of the character $\Theta(\pi)$ of $\pi$ is nonzero (see \cite{HC78} and \cite{MW87}). Since nilpotent orbits $\CO$ of $G_n$ are parametrized by data $(\ul{p}, \ul{q})$, where $\ul{p}$ is partition of $2n$ (or $2n+1$ when $\mathrm{G}_n=\SO_{2n+1}$) and $\ul{q}$ is certain non-degenerate quadratic form (\cite[Section I.6]{Wal01}), we can define another important invariant $\frak{p}(\pi)$ for $\pi$ which is the set of partitions corresponding to $\frak{n}(\pi)$. Under the dominant order of nilpotent orbits or partitions, one can also define $\frak{n}^m(\pi)$ and $\frak{p}^m(\pi)$ to be the maximal elements in $\frak{n}(\pi)$ and $\frak{p}(\pi)$, respectively. We call $\frak{n}^m(\pi)$ the wave front set and $\frak{p}^m(\pi)$ the wave front partitions of $\pi$, respectively.

In general, characterization of the set $\frak{n}^m(\pi)$ is still widely open, though some cases are known. Examples are representations of $\GL_n(F)$, irreducible subquotients for regular principal series of $G_n$ (\cite{MW87}), tempered and anti-tempered representations of pure inner twists of $\SO_{2n+1}(F)$ (\cite{Wal18, Wal20}), certain unramified representations of split connected reductive $p$-adic groups (\cite{Oka21}), and irreducible Iwahori-spherical representations
of split connected reductive $p$-adic groups
with ``real infinitesimal characters" (\cite{CMO21}). It is worth remarking that recently, Tsai (\cite{Tsa24}) constructed an example of representations of $U_7(\mathbb{Q}_3)$ showing that the wavefront set $\frak{p}^m(\pi)$ may not be a singleton.

Given any local Arthur parameter $\psi$ as in \eqref{lap}, the set $\frak{p}^m({\pi}_{\psi})$ can be described as follows.

\begin{thm}\label{wfslinear}
$$\frak{p}^m({\pi}_{\psi})=\{\ul{p}(\psi)^t\}.$$
\end{thm}

Taking the character expansion for the representation $\wt{\pi}_{\psi}$ of the bitorsor $\widetilde{\mathrm{GL}}_N(F)$
at 
\begin{equation}\label{thetaGn}
    \theta_{\widehat{\mathrm{G}}_n}=s_{\widehat{\mathrm{G}}_n} \rtimes \wt{\theta}(N) \in \wt{\GL}_{N}(F),
\end{equation}
(see \cite{Clo87}, also see \cite[Theorem 3.2]{Kon02} and \cite[Theorems 4.20, 4.23]{Var14}), 
where 
$$s_{\widehat{\mathrm{G}}_n}=
\begin{cases}
I_{N}, & \text{ when } \mathrm{G}_n=\SO_{2n+1},\\
\begin{pmatrix}
I_n &&\\
&1&\\
&&-I_n
\end{pmatrix}, & \text{ when } \mathrm{G}_n=\Sp_{2n},\\
\begin{pmatrix}
I_n&\\
&-I_n
\end{pmatrix}, & \text{ when } \mathrm{G}_n= \SO_{2n}^{\alpha},
\end{cases}
$$
we can define the sets $\frak{n}^m(\widetilde{\pi}_{\psi})$ and $\frak{p}^m(\widetilde{\pi}_{\psi})$ 
similarly. 
Note that when $\mathrm{G}_n= \SO_{2n+1}, \SO_{2n}^{\alpha}$, the connected component of the stabilizer of $\theta_{\widehat{\mathrm{G}}_n}$ in $\wt{\mathrm{GL}}_N(F)$ is $\widehat{\mathrm{G}}_n(F)$ and $\frak{n}^m(\widetilde{\pi}_{\psi})$ consists of $F$-rational nilpotent orbits in the Lie algebra of $\widehat{\mathrm{G}}_n(F)$.
When $\RG_n=\Sp_{2n}$,
the connected component of the stabilizer of $\theta_{\widehat{\mathrm{G}}_n}$ in $\wt{\mathrm{GL}}_{2n+1}(F)$ is $\RG_n(F) \times \SO_1$ and $\frak{n}^{m}(\widetilde{\pi}_{\psi})$ consists of $F$-rational nilpotent orbits in the Lie algebra of $\RG_n(F)$.
The choices of $\theta_{\widehat{\mathrm{G}}_n}$ are for the purpose of taking character expansion of certain distributions at $\theta_{\widehat{\mathrm{G}}_n}$ and at its norm.
Then we have the following conjecture regarding the upperbound of the set 
$\frak{p}^m(\widetilde{\pi}_{\psi})$.

\begin{conj}\label{wfsbitorsor}
For any $\ul{p} \in \frak{p}^m(\wt{\pi}_{\psi})$, 
 $$\ul{p}\leq 
 \begin{cases}
(\ul{p}(\psi)^t)_{\widehat{\mathrm{G}}_n}, & \text{ when } \mathrm{G}_n= \SO_{2n+1}, \SO_{2n}^{\alpha},\\
((\ul{p}(\psi)^t)^-)_{\RG_{n}}, & \text{ when } \mathrm{G}_n= \Sp_{2n},\\
 \end{cases}$$
where $(\ul{p}(\psi)^t)_{\widehat{\mathrm{G}}_n}$ is the  $\widehat{\mathrm{G}}_n$-collapse of the partition 
$\ul{p}(\psi)^t$ which is the largest $\widehat{\mathrm{G}}_n$-partition smaller than or equal to $\ul{p}(\psi)^t$, $(\ul{p}(\psi)^t)^-$ is decreasing the smallest part of $\ul{p}(\psi)^t$ by $1$ and  $((\ul{p}(\psi)^t)^-)_{\RG_{n}}$ is the $\RG_n$-collapse of $(\ul{p}(\psi)^t)^-$. 
\end{conj}

We also believe that the following stronger conjecture holds. 

\begin{conj}\label{wfsbitorsor2}
 $$\frak{p}^m(\wt{\pi}_{\psi})=
 \begin{cases}
\{(\ul{p}(\psi)^t)_{\widehat{\mathrm{G}}_n}\}, & \text{ when } \mathrm{G}_n= \SO_{2n+1}, \SO_{2n}^{\alpha},\\
\{((\ul{p}(\psi)^t)^-)_{\RG_{n}}\}, & \text{ when } \mathrm{G}_n= \Sp_{2n}.
 \end{cases}$$
\end{conj}

We also give a generalized version of Conjecture \ref{wfsbitorsor2} (see Conjecture \ref{generalizedconj})  which has its own interests.

\begin{rmk}\label{rmktoconjecture1.3}
(1). By \cite[Theorem 6.3.11]{CM93}, given a local Arthur parameter $\psi$ of $\mathrm{G}_n$, $(\ul{p}(\psi)^t)_{\widehat{\mathrm{G}}_n}$ is always  $\widehat{\mathrm{G}}_n$-special. 

(2). Assume that $\mathrm{G}_n= \SO_{2n+1}, \SO_{2n}^{\alpha}$. When $(\ul{p}(\psi)^t)_{\widehat{\mathrm{G}}_n}=\ul{p}(\psi)^t$, i.e., $\ul{p}(\psi)^t$ is already a $\widehat{\mathrm{G}}_n$-partition, by
Theorem \ref{wfslinear} and 
\cite[Theorem 4.1 (1)]{Kon02}, $\ul{p}(\psi)^t \in \frak{p}^m(\wt{\pi}_{\psi})$.
This covers a large family of local Arthur parameters. Note that 
$\ul{p}(\psi)$ is automatically a $\widehat{\mathrm{G}}_n$-partition, however, $\ul{p}(\psi)^t$ may not be a $\widehat{\mathrm{G}}_n$-partition.
By \cite[Proposition 6.3.7, Theorem 6.3.11]{CM93}, for $\RG_n=\SO_{2n+1}$, $\ul{p}(\psi)^t$ is a $\widehat{\mathrm{G}}_n$-partition if and only if $\ul{p}(\psi)$ is a $\widehat{\mathrm{G}}_n$-special partition; for $\RG_n=\SO^{\alpha}_{2n}$, $\ul{p}(\psi)^t$ is a $\SO_{2n}$-partition if and only if $\ul{p}(\psi)$ is also a $\Sp_{2n}$-partition.

When $(\ul{p}(\psi)^t)_{\widehat{\mathrm{G}}_n}\neq \ul{p}(\psi)^t$, i.e., $\ul{p}(\psi)^t$ is not a $\widehat{\mathrm{G}}_n$-partition, Conjecture \ref{wfsbitorsor2} is expected to be more complicated and is being studied by the authors.

(3). Assume that $\psi$ is tempered, that is for all $1 \leq i \leq r$, $b_i=1$. Then $\ul{p}_{\psi}^t=[N]$.
When $\RG_n=\SO_{2n+1}, \Sp_{2n}$, 
Conjecture \ref{wfsbitorsor2}, hence Conjecture \ref{wfsbitorsor}, has been proved by \cite[Theorem 4.1 (1)]{Kon02} and by \cite[Corollary 6.16]{Var17}. 
When, $\mathrm{G}_n= \SO_{2n}^{\alpha}$, Conjecture \ref{wfsbitorsor2}, hence Conjecture \ref{wfsbitorsor}, has been proved  by \cite[Corollary 6.16]{Var17}.

(4). Assume that 
$\RG_n=\Sp_{2n}$. Then $$((\ul{p}(\psi)^t)^-)_{\RG_{n}}=((\ul{p}(\psi)^-)_{\RG_{n}})^t=\eta_{\frak{so}_{2n+1},\frak{sp}_{2n}}(\ul{p}(\psi)),$$
where $\eta_{\frak{so}_{2n+1},\frak{sp}_{2n}}$ denotes the Barbasch-Vogan duality map from the odd orthogonal partitions to
the symplectic partitions $($see \cite{BV85} and \cite[Lemma 3.3]{Ach03}$)$. 

(5). The choices of $\theta_{\widehat{\mathrm{G}}_n}$ are for the purpose of taking character expansion of certain distributions at $\theta_{\widehat{\mathrm{G}}_n}$ and at its norm.
\end{rmk} 

There is also an interesting question whether all the nilpotent orbits smaller than those in $\frak{n}^m(\pi)$ would occur in $\frak{n}(\pi)$. The answer is yes for the case of $\GL_n$ (\cite{GGS21}), but is unknown in general. 

In this report, we focus on studying the structure of $\frak{p}^m(\pi)$ for irreducible representations $\pi$ of $G_n$ in a local Arthur packet $\widetilde{\Pi}_{\psi}$. It is known that for tempered local Arthur parameters, namely, $b_i=1$ for all $i$, the local Arthur packet is exactly the $L$-packet $\Pi_{\phi_{\psi}}$ corresponding to the $L$-parameter $\phi_{\psi}$. For tempered $L$-packets, the second named author has the following conjecture in general.

\begin{conj}[Shahidi Conjecture]\label{shaconj}
For any quasi-split reductive group $G$, tempered $L$-packets have generic members. 
\end{conj}

This conjecture has been proved for quasi-split classical groups in \cite[Proposition 8.3.2]{Art13}, \cite[Corollary 9.2.4]{Mok15}, based on the global Langlands functoriality \cite{CKPSS04, CPSS11}, \cite{KK04, KK05}, and the automorphic descent \cite{GRS11}; see also the work of \cite{JNQ10}, \cite{JS12}, and \cite{ST15}, via the method of local descent.  Conjecture \ref{shaconj} can be enhanced as follows.

\begin{conj}[Enhanced Shahidi  Conjecture]\label{shaconj2}
For any quasi-split reductive group $G$, local Arthur packets are tempered if and only if they have generic members.  
\end{conj}

We prove the enhanced Shahidi conjecture assuming Conjecture \ref{wfsbitorsor}, see Theorem \ref{mainintro} Part (2) below. We remark that for symplectic and split odd special orthogonal groups, Hazeltine, the first named author and Lo (\cite{HLL22}) have proved Conjecture \ref{shaconj2} without any assumption, using Atobe's refinement on M{\oe}glin's construction of local Arthur packets, which is a different method from that of this paper.

The main goal of this paper is to consider the following conjecture of Jiang which is a natural generalization of the Shahidi conjecture above from tempered local Arthur packets to non-tempered ones, on the characterization of the set $\frak{p}^m(\pi)$ for $\pi$ in local Arthur packets. Note that for a generic representation $\pi$, $\frak{p}^m(\pi)$ contains only regular nilpotent orbits. The global version of this conjecture can be found in \cite[Conjecture 4.2]{Jia14}. We now state the Jiang conjecture as follows. 

\begin{conj}[Jiang Conjecture]\label{cubmfclocal}
Let $\psi$ be a local Arthur parameter of $G_n$ as in \eqref{lap}, and let $\wt{\Pi}_{\psi}$ be the local Arthur packet attached to $\psi$. Then the followings hold.
\begin{enumerate}
\item[(1)] For any partition 
$\ul{p}$ which is not related to $\eta_{{\hat{\frak{g}}_n,\frak{g}_n}}(\ul{p}(\psi))$ and any 
$\pi\in\wt{\Pi}_{\psi}$,
$\ul{p}\notin\frak{p}^m(\pi)$.
\item[(2)] For any partition
$\ul{p}> \eta_{{\hat{\frak{g}}_n,\frak{g}_n}}(\ul{p}(\psi))$ and any 
$\pi\in\wt{\Pi}_{\psi}$,
$\ul{p}\notin\frak{p}^m(\pi)$.
\item[(3)] There exists at least one member $\pi\in\wt{\Pi}_{\psi}$ having the property that
$\eta_{{\hat{\frak{g}}_n,\frak{g}_n}}(\ul{p}(\psi))\in \frak{p}^m(\pi)$.
\end{enumerate}
Here $\eta_{{\hat{\frak{g}}_n,\frak{g}_n}}$ denotes the Barbasch-Vogan duality map from the partitions for the dual group $\widehat{\mathrm{G}}_n(\BC)$ to
the partitions for $G_n$ $($see \cite{BV85} and \cite{Ach03}$)$.
\end{conj}

Our main results are summarized in the following theorem. Recall that given a partition $\ul{p}$ of $M$, where
$$M:= 
\begin{cases}
2n & \text{ when } \mathrm{G}_n=\Sp_{2n}, \SO_{2n}^{\alpha}, \\
2n+1 & \text{ when } \mathrm{G}_n=\SO_{2n+1},
\end{cases}
$$ 
$\ul{p}_{\RG_n}$ denotes the $\RG_n$-collapse of $\ul{p}$, that is, the largest $\RG_n$-partition smaller than or equal to $\ul{p}$, and $\ul{p}^{\RG_n}$ denotes the $\RG_n$-expansion of $\ul{p}$, that is, the smallest special $\RG_n$-partition bigger than or equal to $\ul{p}$. For the definition of special partitions, see \cite[Section 6.3]{CM93}.

\begin{thm}\label{mainintro}
Let $\psi$ be a local Arthur parameter as in \eqref{lap}, with $\ul{p}(\psi) = [b_1^{a_1} b_2^{a_2} \cdots b_r^{a_r}]$ and $b_1 \geq b_2 \geq \cdots \geq b_r$. 
Assume that Conjecture \ref{wfsbitorsor} is true. Then we have the followings. 
\begin{enumerate}
\item Conjecture \ref{cubmfclocal} Part (2) is valid.
\item 
Conjecture \ref{shaconj2} is true.
\item Let $$\underline{p}_1=\left[\big\lfloor \frac{b_1}{2} \big\rfloor^{a_1}
\big\lfloor \frac{b_2}{2} \big\rfloor^{a_2} \cdots \big\lfloor \frac{b_r}{2} \big\rfloor^{a_r}\right]^t,$$
and $n^{*}=\big\lfloor\frac{\sum_{b_i \text{ odd }} a_i}{2}\big\rfloor$. 
    Then Conjecture \ref{cubmfclocal} Part (3) holds for the following cases.
    \begin{enumerate}
        \item When $\mathrm{G}_n=\Sp_{2n}$, and
        \begin{equation}\label{criterion_intro1}
    ([\underline{p}_1\underline{p}_1(2n^*)]^t)_{\Sp_{2n}}=([b_1^{a_1} \cdots b_r^{a_r}]^-)_{\Sp_{2n}}.
\end{equation}
In particular, if 
\begin{enumerate}
    \item $a_r=b_r=1$ and $b_i$ are all even for $1 \leq i \leq r-1$,
    \item or, $b_i$ are all odd,
\end{enumerate}
 then \eqref{criterion_intro1} holds and thus Conjecture \ref{cubmfclocal} Part (3) is valid. 
        \item When $\mathrm{G}_n=\SO_{2n+1}$, and
        \begin{equation}\label{criterion_intro2}
    ([\underline{p}_1\underline{p}_1(2n^*+1)]^t)_{\SO_{2n+1}}=([b_1^{a_1} \cdots b_r^{a_r}]^+)_{\SO_{2n+1}}.
\end{equation}
        In particular, if \begin{enumerate}
            \item $b_1$ is even and $a_1=1$, and $b_i$ are all odd for $2 \leq i \leq r$,
            \item or, $b_i$ are all even, 
        \end{enumerate}
        then \eqref{criterion_intro2} holds  and thus Conjecture \ref{cubmfclocal} Part (3) is valid.
        \item When    $\mathrm{G}_n=\SO_{2n}^{\alpha}$, and
         \begin{equation}\label{criterion_intro3}
    [\underline{p}_1\underline{p}_1(2n^*-1)1]^{\SO_{2n}}=([b_1^{a_1} \cdots b_r^{a_r}]^t)_{\SO_{2n}}.
 \end{equation}
        If all $b_i$ are of the same parity, then \eqref{criterion_intro3} holds and thus Conjecture \ref{cubmfclocal} Part (3) is valid. 
    \end{enumerate}
    
\end{enumerate}
\end{thm}

\begin{rmk}
We should point out that equations \eqref{criterion_intro1}, \eqref{criterion_intro2}, and \eqref{criterion_intro3} of Theorem \ref{mainintro} guarantee the existence of the member in the local Arthur packet suggested by Part (3) of Conjecture \ref{cubmfclocal}. In other words, this part of the conjecture is
equivalent to the validity of these equations for given parameters. 
\end{rmk}

For Parts (1) and (2) of Theorem \ref{mainintro}, we
follow the idea used in \cite[Section 9]{Sha90} which is for the case of endoscopic transfer. More precisely, we apply \cite[Lemmas 9.7, 9.8]{Sha90}, \cite[Lemma 8.5]{Kon02}, and certain dimension identities for nilpotent orbits described in Section 3 below. An ingredient of the proof is a joint result of the authors with Hazeltine and Lo (\cite[Appendix A, Proposition A.4]{LS24}) that for sympletic and special orthogonal groups, all representations in any given local Arthur packet have the same central character. 

For the Part (3) of Theorem \ref{mainintro}, we follow the work of \cite{JS04, Liu11, JL14, JL23} to construct a particular element in each local Arthur packet and study its wave front set, see Section 4 below. 
We remark that when $b_i$ are of mixed parities, conditions \eqref{criterion_intro1} -- \eqref{criterion_intro3} may not always hold, see Remark \ref{mixedpartities} for examples (combining with Theorem \ref{mainthm1}).

Conjecture \ref{cubmfclocal} has analogues for quasi-split unitary groups and non-quasi-split classical groups. 
The method used in this paper is expected to extend to these cases (as long as the local Arthur classification is carried out). 

We remark that Okada (\cite{Oka21}) has proved
Conjecture \ref{cubmfclocal} for certain special local Arthur parameters, namely, those $\psi$ which are trivial on the Weil-Deligne group, and Ciubotaru, Mason-Brown, and Okada (\cite{CMO21}) have computed the wave front set of irreducible Iwahori-spherical representations
of split connected reductive $p$-adic groups
with ``real infinitesimal characters", hence have proved more special cases of Conjecture \ref{cubmfclocal}. 

In the following special cases of Theorem \ref{mainintro}(3), we only need to assume Conjecture \ref{shaconj2}. 

\begin{thm}
Let $\psi$ be a local Arthur parameter of $G_n$ as in \eqref{lap}. Assume that Conjecture \ref{shaconj2} is true. 
    Then Conjecture \ref{cubmfclocal} holds for the following cases.
    \begin{enumerate}
        \item When $\mathrm{G}_n=\Sp_{2n}$,
       $a_1=1$, $b_1=3$, and $b_i=1$ for $2 \leq i \leq r$. 
        \item When $\mathrm{G}_n=\SO_{2n+1}$,
        $a_1=1$, $b_1=2$, and $b_i=1$ for $2 \leq i \leq r$. 
        \item When    $\mathrm{G}_n=\SO_{2n}^{\alpha}$, 
       $a_1=1$, $b_1=3$, and $b_i=1$ for $2 \leq i \leq r$.  
    \end{enumerate}
\end{thm}

\subsection*{Acknowledgements} 
The authors would like to thank Dihua Jiang for his interest and helpful discussions and Bin Xu for helpful communications on Arthur classification of even special orthogonal groups. The authors also would like to thank the referee for helpful comments and suggestions. 

The authors would like to dedicate this article to Steve Kudla on the occasion of his 70th birthday and would like to thank the organizers and the Fields Institute for the invitation. The second author likes to thank Steve for many years of friendships, reflecting upon which leads to many pleasant memories of travels and mathematical discussions. Steve's contributions to the subject have been immense, leading to fundamental progress not only in geometry and number theory but also in automorphic forms and representation theory.

\section{a generalization of Conjecture \ref{wfsbitorsor2}}\label{characlap}

In this section, we describe a generalization of Conjecture \ref{wfsbitorsor2}.

Let $\mathrm{G}$ be a connected reductive group defined over $F$ and let $G=\mathrm{G}(F)$. Let $\theta$ be an automorphism of $G$ of order at most 2. Then $\theta$ is automatically quasi-semisimple (that is, the restriction of $\theta$ to the derived group of $G$ is semisimple) and $G^{\theta}$ is reductive by a theorem of Steinberg (\cite{Ste68}, see also \cite[Theorem 1.1 A]{KS99}). 
Consider the extended group $\widetilde{G}^+=G\rtimes \langle\theta \rangle$, which has two components $G$ and the bitorsor $\widetilde{G}$. Given an irreducible self-dual representation $\pi$ of $G$, that is, $\pi \cong \pi \circ \theta$, let $\widetilde{\pi}^+$ be a unitary extension of $\pi$ to $\widetilde{G}^+$, determined by an intertwining operator $A: \pi \rightarrow \pi \circ \theta$, which is uniquely determined by $\pi$ up to scalars.
Denote by $\widetilde{\pi}$ the restriction of $\widetilde{\pi}^+$ to $\widetilde{G}$.
We refer to \cite[Chapter 2]{LW13} for more details on representations of bitorsors. 
Fix $\wt{\theta}=s
\rtimes \theta$, where $s$ is a semisimple element of $G$. 
Taking the character expansion of $\wt{\pi}$ at $\wt{\theta}$ (see \cite[Theorem 3]{Clo87}, also see \cite[Theorem 3.2]{Kon02} and \cite[Theorems 4.20, 4.23]{Var14}), we define the wave front set $\frak{n}^m(\widetilde{\pi})$, which consists of $G^{\wt{\theta}}$-nilpotent orbits in its Lie algebra. 
Then we have the following conjecture about the relation between $\frak{n}^m(\pi)$ and $\frak{n}^m(\widetilde{\pi})$.  
 
\begin{conj}\label{generalizedconj}
    $$\frak{n}^m(\widetilde{\pi}) = \left\{ (\CO_{G^{\wt{\theta}}})^{G^{\wt{\theta}}} | \CO \in \frak{n}^m(\pi)\right\},$$
    where $\CO_{G^{\wt{\theta}}}$ is the  $G^{\wt{\theta}}$-collapse of $\CO$ which is the largest $G^{\wt{\theta}}$-nilpotent orbit smaller than or equal to $\CO$, and $(\CO_{G^{\wt{\theta}}})^{G^{\wt{\theta}}}$ is the
    $G^{\wt{\theta}}$-expansion of $\CO_{G^{\wt{\theta}}}$ which is the smallest $G^{\wt{\theta}}$-special nilpotent orbit bigger than or equal to $\CO_{G^{\wt{\theta}}}$. Here, special nilpotent orbits correspond to special representations of Weyl groups via the Springer correspondence.
\end{conj}

When $\theta$ is trivial, we are reduced to the connected group case, and by considering the character expansion at another semisimple element (\cite{HC78}), the conjecture still seems to be interesting and open. 

To obtain Conjecture \ref{wfsbitorsor2} which is a special case of Conjecture \ref{generalizedconj}, we take the character expansion for the representation $\wt{\pi}_{\psi}$ of the bitorsor $\widetilde{\mathrm{GL}}_N(F)$ at the point $\theta_{\widehat{\mathrm{G}}_n}$ as in \eqref{thetaGn}. Note that when $\mathrm{G}_n= \SO_{2n+1}, \SO_{2n}^{\alpha}$, the connected component of the stabilizer of $\theta_{\widehat{\mathrm{G}}_n}$ in $\mathrm{GL}_N$ is $\widehat{\mathrm{G}}_n$, and by \cite[Theorem 6.3.11]{CM93}, given a local Arthur parameter $\psi$ of $\mathrm{G}_n$, $(\ul{p}(\psi)^t)_{\widehat{\mathrm{G}}_n}$ is always $\widehat{\mathrm{G}}_n$-special. Here, we are considering partitions, not rational nilpotent orbits, and hence, abusing the notation, we write $(\ul{p}(\psi)^t)_{\widehat{\mathrm{G}}_n}$ for convenience, instead of $(\ul{p}(\psi)^t)_{\mathrm{GL}_N^{\theta_{\widehat{\mathrm{G}}_n}}}$. 
When $\RG_n=\Sp_{2n}$,
the connected component of the stabilizer of $\theta_{\widehat{\mathrm{G}}_n}$ in $\wt{\mathrm{GL}}_{2n+1}(F)$ is $\RG_n(F) \times \SO_1$ and $\frak{n}^{m}(\widetilde{\pi}_{\psi})$ consists of $F$-rational nilpotent orbits in the Lie algebra of $\RG_n(F)$. In this case, one can see that $(\ul{p}(\psi)^t)_{\mathrm{GL}_N^{\theta_{\widehat{\mathrm{G}}_n}}}$ is exactly $((\ul{p}(\psi)^t)^-)_{\RG_{n}}$, which is equal to $\eta_{\frak{so}_{2n+1},\frak{sp}_{2n}}(\ul{p}(\psi))$ by Remark \ref{rmktoconjecture1.3}(4), hence it is special.


\section{Dimension identities for nilpotent orbits}\label{dimidentities}

In this section, we state certain dimension identities for nilpotent orbits which are important ingredients for our main results.

Recall that $\psi = \bigoplus_{i=1}^r \phi_i \otimes S_{m_i} \otimes S_{n_i}$, where $\phi_i$ is of dimension $k_i$, $a_i=k_im_i$, $b_i=n_i$, and $2n+1=\sum_{i=1}^r a_ib_i$. 
For $1 \leq i \leq r$, let $\psi_i=\phi_i \otimes S_{m_i} \otimes S_{n_i}$. Let $\{1,2,\ldots, r\}=I \dot\cup J$ such that $I=\{i|s_{\psi_i}=1\}$, $J=\{i|s_{\psi_i}=-1\}$.
By definition, $s_{\psi}=\psi\left(1,1, \begin{pmatrix}
-1 &0 \\
0&-1
\end{pmatrix}\right)$, and
\begin{align*}
    &\,\phi_i \otimes S_{m_i} \otimes S_{b_i}
    \left(1,1, \begin{pmatrix}
-1 &0 \\
0&-1
\end{pmatrix}\right)\\
=&\,\begin{pmatrix}
(-1)^{b_i-1}I_a&&&\\
&(-1)^{b_i-3}I_a&&\\
&&\ddots&\\
&&&(-1)^{1-b_i}I_a
\end{pmatrix}.
\end{align*}
Hence, 
$I=\{i|b_i \text{ odd}\}$, $J=\{i|b_i \text{ even}\}$, and 
automatically, we have $\sum_{i\in I}a_ib_i=N_{n_1}$, $\sum_{j \in J}a_jb_j=2n_2$, and $n=n_1+n_2$. 
Moreover, $s_{\psi}$ has the form $$\begin{pmatrix}
-I_{n_2}&&\\
&I_{N_{n_1}}&\\
&&-I_{n_2}
\end{pmatrix}.$$

We now apply \cite[Theorem 2.2.1, Part (b)]{Art13} with $s=s_{\psi}=x$. It is clear that the stabilizer of $s$ in $\widehat{\mathrm{G}}_n(\BC)$ is
\begin{align*}
    \SO_{2n_1+1}(\BC) \times \SO_{2n_2}(\BC) & \text{ when } \mathrm{G}_n=\Sp_{2n},\\
    \Sp_{2n_1}(\BC) \times \Sp_{2n_2}(\BC) & \text{ when } \mathrm{G}_n=\SO_{2n+1},\\
    \SO_{2n_1}(\BC) \times \SO_{2n_2}(\BC) & \text{ when } \mathrm{G}_n=\SO_{2n}^{\alpha},
\end{align*}
and the corresponding endoscopic group $\mathrm{G}'=\mathrm{G}_1' \times \mathrm{G}_2'$ of $\mathrm{G}_n$ is
\begin{align*}
    \Sp_{2n_1} \times \SO_{2n_2}^{\beta} & \text{ when } \mathrm{G}_n=\Sp_{2n},\\
    \SO_{2n_1+1} \times \SO_{2n_2+1} & \text{ when } \mathrm{G}_n=\SO_{2n+1},\\
    \SO_{2n_1}^{\gamma_1} \times \SO_{2n_2}^{\gamma_2} & \text{ when } \mathrm{G}_n=\SO_{2n}^{\alpha},
\end{align*}
where $\beta, \gamma_1, \gamma_2$ are square classes in $F$, and $\gamma_1\gamma_2=\alpha$. 
By \cite[Theorem 2.2.1]{Art13}, we have the following distribution identity
$$\sum_{\pi \in \wt{\Pi}_{\psi}} \langle s_{\psi} x, \pi \rangle f_G(\pi) = f'(\psi')=\tr (\wt{\pi}_{\psi^1}(\wt{f}^1)) \tr (\wt{\pi}_{\psi^2}(\wt{f}^2)),$$
where $f \in \wt{\CH}(G_n)$, $f' \in \wt{\CS}(G')$ is the transfer of $f$ to $G'$ with the assumption that $f'=f'^1 \otimes f'^2$, $f'^i \in \wt{\CS}(G'_i)$, $i=1, 2$.
Let 
$$\wt{f}^1 \in \begin{cases}
\wt{\CH}(2n_1+1) & \text{ when } \mathrm{G}_n=\Sp_{2n},\\
\wt{\CH}(2n_1) & \text{ when } \mathrm{G}_n=\SO_{2n+1}, \SO_{2n}^{\alpha},
\end{cases}
$$
and 
$\wt{f}^2 \in \wt{\CH}(2n_2)$, transferring to $f'^1$ and $f'^2$, via surjective maps in \cite[Corollary 2.1.2]{Art13}
$$\iota_{G_k}: \wt{\CH}(N) \rightarrow \wt{\CS}(G_k), k=n_1, n_2,$$
respectively. Here, $\psi'$ is the factor through of $\psi$ to ${}^L\mathrm{G}'$, and $\psi^1=\sum_{i\in I}\psi_i$, $\psi^2 = \sum_{j \in J}\psi_j$. 

Take the character expansions as in \cite{HC78} for connected groups at
$$c_{G_n}= 
\begin{cases}
I_{2n+1} & \text{ when } \mathrm{G}_n=\SO_{2n+1},\\
-I_{2n} & \text{ when } \mathrm{G}_n=\Sp_{2n}, \SO_{2n}^{\alpha}.
\end{cases}
$$ 
Take the character expansions as in \cite{Clo87} for disconnected groups (at $\theta_{\widehat{\mathrm{G}}_i'}$, see \eqref{thetaGn}, $i=1,2$). Then, we have the following equality: 
\begin{align}\label{sec6equ1sp2n}
\begin{split}
    &\,  \sum_{\pi \in \wt\Pi_{\psi}} \langle 1, \pi \rangle \sum_{\CO \in \CN_{\frak{g}_{n}}} \omega_{\pi}(c_{\RG_n})
c_{\CO}(\pi)\hat{\mu}_\CO(R(c_{\RG_n})f) \\
= & \, \sum_{\CO \in \CN_{\frak{g}_{n}}}
c_{\CO}(\wt{\Pi}_{\psi})\hat{\mu}_\CO(R(c_{\RG_n})f)\\
= &\,  \begin{cases}
    \left(\sum_{\CO \in \CN_{{\frak{g}}_1'}}
c_{\CO}(\wt{\pi}_{\psi^1})\hat{\mu}_\CO({{}\wt{f}^1_{\theta_{\widehat{\mathrm{G}}_1'}}})\right) \left(\sum_{\CO \in \CN_{\hat{\frak{g}}_2'}}
c_{\CO}(\wt{\pi}_{\psi^2})\hat{\mu}_\CO({{}\wt{f}^2_{\theta_{\widehat{\mathrm{G}}_2'}}})\right), \\
\text{ when } \mathrm{G}_n=\Sp_{2n},\\
\left(\sum_{\CO \in \CN_{\hat{\frak{g}}_1'}}
c_{\CO}(\wt{\pi}_{\psi^1})\hat{\mu}_\CO({{}\wt{f}^1_{\theta_{\widehat{\mathrm{G}}_1'}}})\right) \left(\sum_{\CO \in \CN_{\hat{\frak{g}}_2'}}
c_{\CO}(\wt{\pi}_{\psi^2})\hat{\mu}_\CO({{}\wt{f}^2_{\theta_{\widehat{\mathrm{G}}_2'}}})\right), \\
\text{ when } \mathrm{G}_n=\SO_{2n+1}, \SO_{2n}^{\alpha},\\
\end{cases}
\end{split}
\end{align}
where $\CN_{\frak{g}}$ denotes the set of $F$-rational nilpotent orbits in a Lie algebra $\frak{g}$, {{}$\wt{f}_{\theta_{\widehat{\mathrm{G}}_i'}}^i$ is the Harish-Chandra descent of $\widetilde{f}^i$ (see for example \cite[Section 3.1]{Kon02}), $i=1,2$,} 
and $c_{\CO}(\wt{\Pi}_{\psi})=\sum_{\pi \in \wt\Pi_{\psi}} \langle 1, \pi \rangle 
c_{\CO}(\pi)$. Note that by Conjecture \ref{wfsbitorsor} and Remark \ref{rmktoconjecture1.3}(4), 
the maximal nilpotent orbits that occur on the right-hand side are those corresponding to partitions $([\prod_{i\in I}b_i^{a_i}]^t)_{\widehat{\mathrm{G}}_1'}$ and $([\prod_{j\in J}b_j^{a_j}]^t)_{\widehat{\mathrm{G}}_2'}$, or $\eta_{{\hat{\frak{g}}_1',\frak{g}_1'}}(\ul{p}(\psi^1))$ if $\G_1'=\Sp_{2n_1}$.  

We have the following dimension identities for nilpotent orbits, which are ingredients in the proof of Theorem \ref{mainintro}(1) and (2).

\begin{lem}\label{keylemma2}
We have the following dimension identity:
\begin{align}\label{equ1:keylemma2}
\begin{split}
    &\, \dim(\frak{g}_{n}) - \dim_{\frak{g}_{n}}(\eta_{\hat{\frak{g}}_{n}, \frak{g}_{n}}([b_1^{a_1} \cdots b_r^{a_r}])) \\
    = &\, \dim(\hat{\frak{g}}_{1}') + \dim(\hat{\frak{g}}_{2}') -\dim_{\hat{\frak{g}}_{1}'}(([\prod_{i\in I}b_i^{a_i}]^t)_{\widehat{\mathrm{G}}_1'})-\dim_{\hat{\frak{g}}_{2}'}(([\prod_{j\in J}b_j^{a_j}]^t)_{\widehat{\mathrm{G}}_2'}).
    \end{split}
\end{align}
\end{lem}

\begin{lem}\label{keylemma1}
We have the following dimension identity:
\begin{equation}\label{equ:keylemma}
    \dim_{\frak{g}_{n}}(\eta_{\hat{\frak{g}}_{n}, \frak{g}_{n}}([b_1^{a_1} \cdots b_r^{a_r}])) = \dim_{\hat{\frak{g}}_{n}}(([b_1^{a_1} \cdots b_r^{a_r}]^t)_{\widehat{\mathrm{G}}_n}),
\end{equation}
where $([b_1^{a_1} \cdots b_r^{a_r}]^t)_{\widehat{\mathrm{G}}_n}$ is the $\widehat{\mathrm{G}}_n$-collapse of the partition $[b_1^{a_1} \cdots b_r^{a_r}]^t$.
\end{lem}

\section{Construction of elements in local Arthur packets}\label{nonvan}

In this section, we construct a particular element in each local Arthur packet which plays an important role towards proving Part (3) of Conjecture \ref{cubmfclocal}. Then we state certain results on its wave front set. 

Given a local Arthur parameter $$\psi: W_F \times \SL_2(\mathbb{C}) \times \SL_2(\mathbb{C}) \rightarrow {}^L\mathrm{G}_{n}$$
$$\psi = \bigoplus_{i=1}^r \phi_i \otimes S_{m_i} \otimes S_{n_i},$$
as in \eqref{lap}, we have
$$\phi_{\psi}=\bigoplus_{i=1}^r 
\bigoplus_{j=-\frac{n_i-1}{2}}^{\frac{n_i-1}{2}} |w|^j\phi_i(w) \otimes S_{m_i} (x).$$

Let 
$$\phi^{(t)}_{\psi}=
\bigoplus_{i=1, n_i \text{ odd }}^r 
\phi_i(w) \otimes S_{m_i} (x),$$
and 
$$\phi^{(n)}_{\psi}=\bigoplus_{i=1}^r 
\bigoplus_{j=-\frac{n_i-1}{2}, j\neq 0}^{\frac{n_i-1}{2}} |w|^j\phi_i(w) \otimes S_{m_i} (x).$$
Then $\phi^{(t)}_{\psi} \in \Phi^{(t)}(G_{n^{*}})$ and 
$\phi^{(n)}_{\psi} \in \Phi(G'_{n-n^{*}})$, where 
$$n^{*}=\bigg\lfloor\frac{\sum_{n_i \text{ odd }} k_im_i}{2}\bigg\rfloor.$$  
By \cite{JS04, Liu11, JL14, JL23}, there exists
$\sigma^{(t)} \in \Pi^{(tg)}(G_{n^{*}})$, where ``$tg$" means tempered generic, such that
\begin{equation} \label{equ692}
\iota(\sigma^{(t)}) = \phi^{(t)}_{\psi}. 
\end{equation}
Using the local Langlands reciprocity map $r$ for general linear groups, define
\begin{equation} \label{equ702}
\Sigma_i^j = [v^{j-\frac{m_i-1}{2}}r(\phi_i), v^{j+\frac{m_i-1}{2}}r(\phi_i)], 1 \leq i
\leq r, 0 < j \leq \frac{n_i-1}{2}.
\end{equation}
Shuffle the set $\{0 <  j \leq \frac{n_i-1}{2} | 1 \leq i \leq r\}$ as $\{j_1, \ldots, j_{n'}\}$ such that $j_1 \geq j_2 \geq \cdots \geq j_{n'}$,  where $n'=\sum_{i=1}^r \frac{n_i-1}{2}$. 
Let $\sigma$ be the Langlands quotient
of the induced representation
$$\delta(\Sigma_1) \times \delta(\Sigma_2) \times \cdots \times 
\delta(\Sigma_{n'}) \rtimes \sigma^{(t)}.$$ 
Then $\sigma$ is an irreducible representation with local Langlands parameter $\phi_{\psi}$ and $\sigma$ is in the local Arthur packet corresponding to $\psi$. 

Recall that for $1 \leq i \leq r$, we let $a_i=k_im_i$ and $b_i=n_i$. Let $\underline{p}_1=[\lfloor \frac{b_1}{2} \rfloor^{a_1}
\lfloor \frac{b_2}{2} \rfloor^{a_2} \cdots \lfloor \frac{b_r}{2} \rfloor^{a_r}]^t$. Recall that $n^{*}=\big\lfloor\frac{\sum_{b_i \text{ odd }} a_i}{2}\big\rfloor.$ Then we have the following result on the set $\frak{p}(\sigma)$, which is proved by applying results in \cite{JLS16, GGS17}.

\begin{lem}\label{raise}
\begin{enumerate}
    \item $[\underline{p}_1\underline{p}_1(2n^*)]^{\mathrm{G}_n} \in \frak{p}(\sigma), \text{ when } \mathrm{G}_n=\Sp_{2n}$,\\
\item $[\underline{p}_1\underline{p}_1(2n^*+1)]^{\mathrm{G}_n} \in \frak{p}(\sigma), \text{ when } \mathrm{G}_n=\SO_{2n+1}$,\\
    \item $[\underline{p}_1\underline{p}_1(2n^*-1)1]^{\mathrm{G}_n} \in \frak{p}(\sigma), \text{ when } \mathrm{G}_n= \SO_{2n}^{\alpha}$.
\end{enumerate}
\end{lem}

Combining with Lemma \ref{raise}, the following theorem proves certain cases for Part (3) of Conjecture \ref{cubmfclocal}. 

\begin{thm}\label{mainthm1}
Let $\psi$ be a local Arthur parameter as in \eqref{lap}, with $\ul{p}(\psi) = [b_1^{a_1} b_2^{a_2} \cdots b_r^{a_r}]$ and $b_1 \geq b_2 \geq \cdots \geq b_r$. 
\begin{enumerate}
    \item When $\mathrm{G}_n=\Sp_{2n}$,
    \begin{equation}\label{oriequ}
    [\underline{p}_1\underline{p}_1(2n^*)]^{\Sp_{2n}}=\eta_{\frak{so}_{2n+1}, \frak{sp}_{2n}}([b_1^{a_1} \cdots b_r^{a_r}])
    \end{equation}
if and only if 
\begin{equation}\label{criterion}
    ([\underline{p}_1\underline{p}_1(2n^*)]^t)_{\Sp_{2n}}=([b_1^{a_1} \cdots b_r^{a_r}]^-)_{\Sp_{2n}},
\end{equation}
where $[b_1^{a_1} \cdots b_r^{a_r}]^-=[b_1^{a_1} \cdots b_r^{a_r-1}(b_r-1)]$. 
In particular, if 
\begin{enumerate}
    \item[(i)] $a_r=b_r=1$ and $b_i$ are all even for $1 \leq i \leq r-1$,
    \item[(ii)] or, $b_i$ are all odd,
\end{enumerate}
 then \eqref{criterion} holds. 

 \item When $\mathrm{G}_n=\SO_{2n+1}$, 
 \begin{equation}\label{oriequ2}
     [\underline{p}_1\underline{p}_1(2n^*+1)]^{\SO_{2n+1}}=\eta_{\frak{sp}_{2n}, \frak{so}_{2n+1}}([b_1^{a_1} \cdots b_r^{a_r}])
 \end{equation}
if and only if 
\begin{equation}\label{criterion2}
    ([\underline{p}_1\underline{p}_1(2n^*+1)]^t)_{\SO_{2n+1}}=([b_1^{a_1} \cdots b_r^{a_r}]^+)_{\SO_{2n+1}},
\end{equation}
where $[b_1^{a_1} \cdots b_r^{a_r}]^+=[(b_1+1)b_1^{a_1-1} \cdots b_r^{a_r}]$. 
        In particular, if \begin{enumerate}
    \item[(i)] $b_1$ is even and $a_1=1$, and $b_i$ are all odd for $2 \leq i \leq r$,
    \item[(ii)] or, $b_i$ are all even,
\end{enumerate}
then 
\eqref{criterion2} holds. 

 \item Assume $\mathrm{G}_n=\SO_{2n}^{\alpha}$.
 If all $b_i$ are of the same parity, then 
 \begin{equation}\label{criterion3}
    [\underline{p}_1\underline{p}_1(2n^*-1)1]^{\SO_{2n}}=\eta_{\frak{o}_{2n}, \frak{o}_{2n}}([b_1^{a_1} \cdots b_r^{a_r}]).
 \end{equation}
\end{enumerate}
\end{thm}

We remark that the identities \eqref{criterion} and \eqref{criterion2} are relatively easier to check than \eqref{oriequ} and \eqref{oriequ2}.

\begin{rmk}\label{mixedpartities}
When $b_i$ are of mixed parities, conditions \eqref{oriequ}, \eqref{oriequ2}, and \eqref{criterion3} may not always hold.  We give some examples as follows. 

When $\mathrm{G}_n=\Sp_{2n}$, if $\ul{p}(\psi)=[3^32^2]$, then 
$$[\underline{p}_1\underline{p}_1(2n^*)]^{\mathrm{G}_n}=\eta_{\hat{\frak{g}}_{n}, \frak{g}_{n}}([3^32^2])=[5^22];$$
if $\ul{p}(\psi)=[2^21^3]$, then 
$$[\underline{p}_1\underline{p}_1(2n^*)]^{\mathrm{G}_n}=[2^3]<\eta_{\hat{\frak{g}}_{n}, \frak{g}_{n}}([2^21^3])=[42].$$

When $\mathrm{G}_n=\SO_{2n+1}$, if
$\ul{p}(\psi)=[3^22^3]$, then
$$[\underline{p}_1\underline{p}_1(2n^*+1)]^{\mathrm{G}_n}=\eta_{\hat{\frak{g}}_{n}, \frak{g}_{n}}([3^22^3])=[5^23];$$
if $\ul{p}(\psi)=[2^31^2]$, then
$$[\underline{p}_1\underline{p}_1(2n^*+1)]^{\mathrm{G}_n}=[3^3]<\eta_{\hat{\frak{g}}_{n}, \frak{g}_{n}}([2^31^2])=[531].$$

When $\mathrm{G}_n=\SO_{2n}^{\alpha}$, if $\ul{p}(\psi)=[3^32^21]$, then 
$$[\underline{p}_1\underline{p}_1(2n^*-1)1]^{\mathrm{G}_n}=\eta_{\hat{\frak{g}}_{n}, \frak{g}_{n}}([3^32^2])=[5^231];$$
if $\ul{p}(\psi)=[32^21^3]$, then 
$$[\underline{p}_1\underline{p}_1(2n^*-1)1]^{\mathrm{G}_n}=[3^31]<\eta_{\hat{\frak{g}}_{n}, \frak{g}_{n}}([2^21^3])=[631].$$
\end{rmk}

In general, we have the following inequality of partitions which is consistent with Conjecture \ref{cubmfclocal}.

\begin{prop}\label{propinequality}
Let $\mathrm{G}_n=\Sp_{2n}, \SO_{2n+1}$, and let $\ul{p}=[\underline{p}_1\underline{p}_1(2n^*)]^{\mathrm{G}_n}$, or $[\underline{p}_1\underline{p}_1(2n^*+1)]^{\mathrm{G}_n}$, then
$$\ul{p}\leq \eta_{\hat{\frak{g}}_{n}, \frak{g}_{n}}([b_1^{a_1} \cdots b_r^{a_r}]).$$ 
\end{prop}

\begin{rmk}
When $\mathrm{G}_n=\SO_{2n}^{\alpha}$, we also expect that 
$$[\underline{p}_1\underline{p}_1(2n^*-1)1]^{\SO_{2n}}\leq
\eta_{\frak{so}_{2n}, \frak{so}_{2n}}([b_1^{a_1} \cdots b_r^{a_r}])=([b_1^{a_1} \cdots b_r^{a_r}]^t)_{\SO_{2n}}.$$
It is an interesting question to completely describe the conditions on $a_i, b_i$ such that the equalities \eqref{criterion}, \eqref{criterion2}, and \eqref{criterion3} hold. 
\end{rmk}


\end{document}